\documentclass[11pt]{article}
\usepackage{bridges}
\usepackage{amsfonts,amssymb,amsthm,eucal,amsmath}
\usepackage{graphicx}
\usepackage[T1]{fontenc}
\usepackage{latexsym,url}
\usepackage{enumerate}

\usepackage{pinlabel} 
\usepackage{subfig}
\usepackage{wrapfig}

\theoremstyle{definition}

\theoremstyle{remark}

\def\H{\mathbb{H}}
\def\R{\mathbb{R}}
\newcommand{\BB}{\mathbb{B}}
\newcommand{\TT}{\mathbb{T}}

\title{Sculptures in $S^3$} 

\author{
\begin{tabular}{cc}
Saul Schleimer & Henry Segerman\footnote{This work is in the
public domain.}\\
Mathematics Institute & Department of Mathematics and Statistics\\
University of Warwick & University of Melbourne\\
Coventry CV4 7AL & Parkville VIC 3010\\
United Kingdom & Australia\\
s.schleimer@warwick.ac.uk & segerman@unimelb.edu.au
\end{tabular}}

\date{}

\begin{document}
\maketitle
\begin{abstract}
We construct a number of sculptures, each based on a geometric design
native to the three-dimensional sphere.  Using stereographic
projection we transfer the design from the three-sphere to ordinary
Euclidean space.  All of the sculptures are then fabricated by the 3D
printing service Shapeways.
\end{abstract}

\section{Introduction}

\begin{figure}[htb]
\centering
\subfloat[Stereographic projection from $S^1 - \{N\}$ to $\R^1$. ]{
\labellist
\large\hair 2pt
\pinlabel $N$ at 135 213
\pinlabel $\frac{x}{1-y}$ at 102 86
\pinlabel $\mathbb{R}^1$ at -8 107
\normalsize
\pinlabel $(x,y)$ at 46 18
\endlabellist
\includegraphics[width=0.30\textwidth]{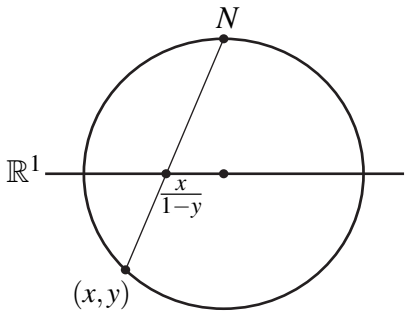}
\label{stereo_proj_diagram}}
\qquad
\subfloat[Two-dimensional stereographic projection applied to the
Earth. Notice that features near the north pole are very large in the
image.]{
\includegraphics[width=0.40\textwidth]{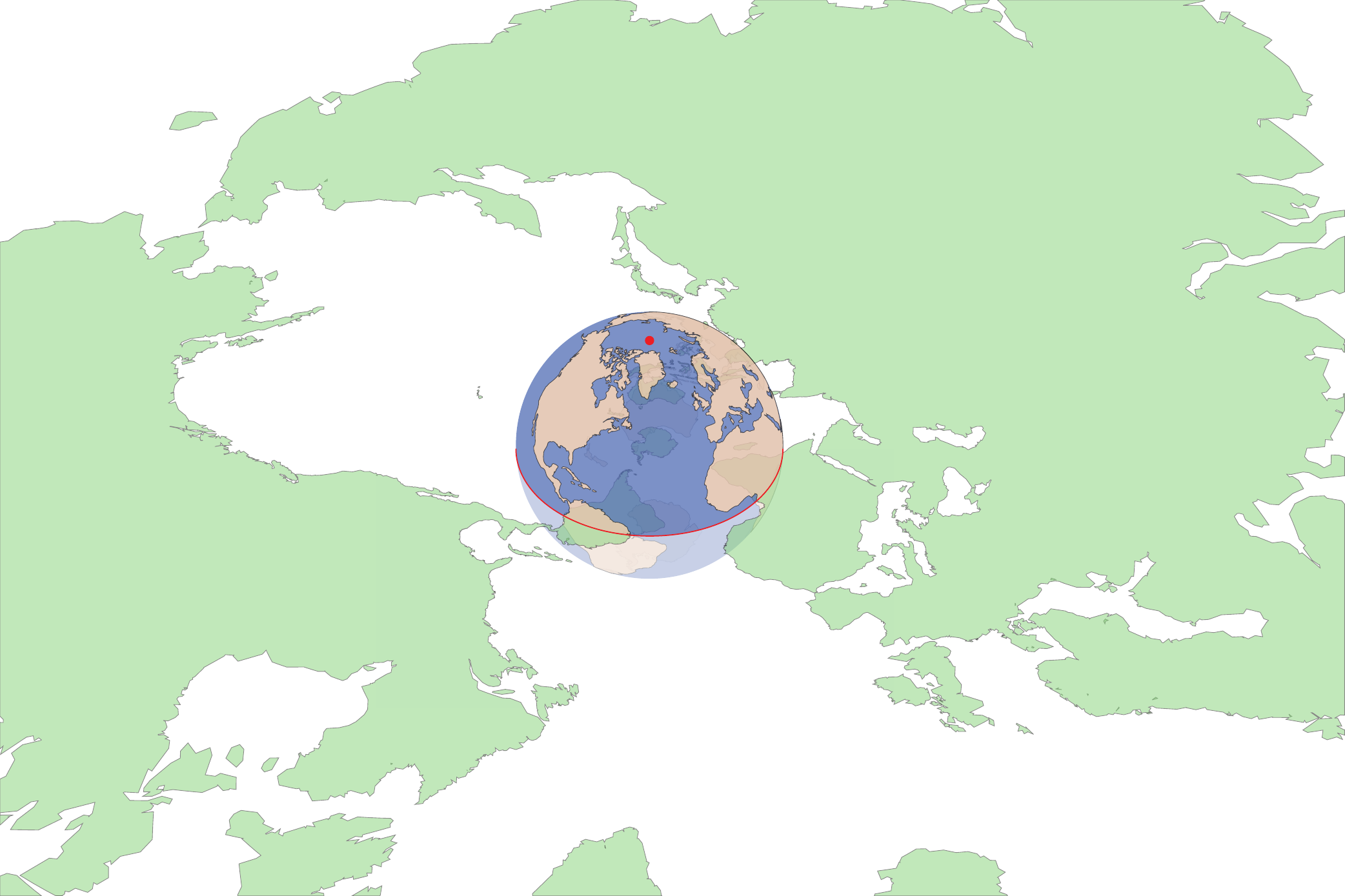}
\label{stereo_proj_earth.pdf}}
\caption{Stereographic projection.}
\label{Stereographic projection}
\end{figure}

The three-sphere, denoted $S^3$, is a three-dimensional analog of the
ordinary two-dimensional sphere, $S^2$.  In general, the
$n$--dimensional sphere is a subset of Euclidean space, $\R^{n+1}$, as
follows:
\[
S^n = \{ (x_0, x_1, \ldots,x_n) \in \R^{n+1} \mid 
                        x_0^2+x_1^2+\cdots+x_n^2=1 \}. 
\]
Thus $S^2$ can be seen as the usual unit sphere in $\R^3$.
Visualising objects in dimensions higher than three is non-trivial.
However for $S^3$ we can use stereographic projection to reduce the
dimension from four to three.  Let $N=(0,\ldots,0,1)$ be the north
pole of $S^n$.  We define {\bf stereographic projection} $\rho : S^n -
\{N\}\to\R^n$ by
\[
\rho(x_0, x_1, \ldots,x_n) = \left(\frac{x_0}{1-x_n},
    \frac{x_1}{1-x_n}, \ldots,\frac{x_{n-1}}{1-x_n}\right).
\]
See~\cite[page 27]{Beardon1983}.  Figure~\ref{stereo_proj_diagram}
displays stereographic projection in dimension one.  For any point
$(x,y) \in S^1 - \{ N \}$ draw the straight line $L$ through $N$ and
$(x,y)$.  Then $L$ meets $\R^1$ at a single point; this is
$\rho(x,y)$.  Notice that the figure is also a two-dimensional
cross-section of stereographic projection in any dimension.
Additionally, there is nothing special about the choice of $N =
(0,\ldots,0,1)$.  We may alter the formula so that any point in $S^3$
becomes the \textbf{projection point}.

By adding in a point at infinity corresponding to the north pole,
stereographic projection extends to a homeomorphism from $S^n$ to $\R^n
\cup \{\infty\}$.  So we may use stereographic projection to
represent, in $\R^3$, objects that live in $S^3$.

\section{The geometry of $S^3$}

A generic plane in $\R^4$, meeting $S^3$, meets $S^3$ in a circle.
The following \textbf{circline property} is fundamental: stereographic
projection maps any circle $C \subset S^3$ to a circle or line in
$\R^3$.  Accordingly we use the word \textbf{circline} as a shorthand
for circles and lines in $\R^3$.  See~\cite[Section~3.2]{Beardon1983}
for a more general discussion.
Note that a circle $C$ of $S^3$ maps to a line in $\R^3$ if and only
if $C$ meets the projection point.

Any plane, meeting the origin in $\R^4$, cuts $S^3$ in a \textbf{great
circle}.  The great circles are the geodesics, or locally shortest
paths, in the geometry on $S^3$.  Just as for the usual sphere, $S^2$,
two distinct great circles meet at two points: say at $x \in \R^4$ and
also at the antipodal point $-x$.

Stereographic projection is \textbf{conformal}: if two circles in
$S^3$ intersect at a given angle then the corresponding circlines in
$\R^3$ meet at the same angle.  So stereographic projection preserves
angles~\cite[Section~3.2]{Beardon1983}.  Note that lengths are not
preserved; as shown in Figure~\ref{stereo_proj_earth.pdf} the
distortion of length becomes infinite as we approach the projection
point.  However, this defect is unavoidable; there is no isometric
embedding of any open subset of the three-sphere into $\R^3$.

\begin{wrapfigure}[16]{r}{0.5\textwidth}

\centering
\labellist
\hair 2pt
\pinlabel $-j$ at 129 181
\pinlabel $i$ at 174 121
\pinlabel $1$ at 221 163
\pinlabel $k$ at 221 253
\pinlabel $-k$ at 245 47
\pinlabel $-i$ at 254 194
\pinlabel $j$ at 335 140
\endlabellist
\includegraphics[width=0.5\textwidth]{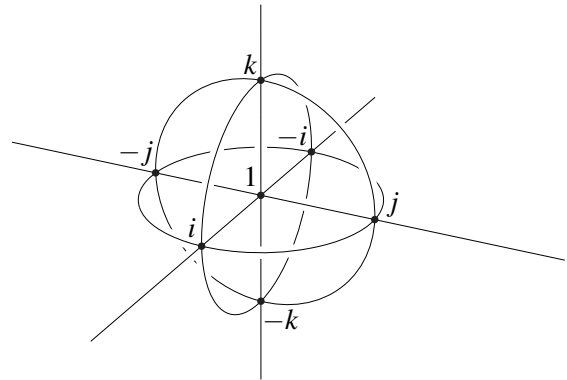}
\caption{The unit quaternions in $S^3$ stereographically projected to
$\R^3$ from the projection point $-1$.}
\label{axes.pdf}
\end{wrapfigure}

\paragraph{The quaternionic picture of $S^3$}
In order to get a sense of the shape of $S^3$, it is useful to have
some landmarks.  A good way to do this is to view $S^3$ in terms of
the unit quaternions~\cite{ConwaySmith2003}.  The quaternions are an
extension of the complex numbers, from two dimensions to four.  A
quaternion is a formal sum $a + bi + cj + dk$ where $a,b,c,d\in\R$ and
where $i,j,k$ are non-commuting symbols satisfying
\[
i^2 = j^2 = k^2 = ijk = -1.
\]
The set of quaternions is called $\H$ in honour of Hamilton, its
discoverer.  There is a natural bijection between $\R^4$ and
$\mathbb{H}$ via $(a,b,c,d) \mapsto a + bi+cj+dk$.  So we may view
$S^3$ as the set of \textbf{unit quaternions}: those with length
$|a+bi+cj+dk| = \sqrt{a^2+b^2+c^2+d^2}$ equal to one.  Once this is
established the points $\pm 1, \pm i, \pm j, \pm k$ serve as our
landmarks.  See Figure~\ref{axes.pdf}.  All of the circlines shown
correspond to great circles in $S^3$ with particularly nice
quaternionic expressions.

\paragraph{The isometries of $S^3$}

The isometries of $S^1$ are the rigid motions of $\R^2$ that fix the
origin, namely rotations and reflections.  Under composition, these
form the \textbf{orthogonal group} $O(2)$.
Analogously, the isometries of $S^n$ form the group $O(n)$.  The unit
quaternions can be realised as a subgroup of $O(4)$ in the following
manner.  As above we identify $\H$ and $\R^4$.  For $q \in \H$ with
$|q|=1$, the map $f_q:\H\to\H$ given by $f_q(x)=q\cdot x$ is an
element of $O(4)$.  So, if we want to move the point $a$ to the point
$b$ in $S^3$, then one way to achieve this is to apply the isometry
corresponding to the quaternion $b\cdot a^{-1}$.

An application of this technique is to adjust the stereographic
projection of a subset of $S^3$.  If $F \subset S^3$ is a surface
then, as $q$ varies, the image of $q \cdot F$ in $\R^3$ changes
dramatically.  Equivalently we can think of this as moving the
projection point.  

\section{Designs in $S^3$}

\subsection{Four-dimensional polytopes}
\label{Sec:Polytopes}

Suppose that $\sigma \subset \R^n$ is a finite set.  Then $P =
P(\sigma)$, the convex hull of $\sigma$, is a
\textbf{polytope}~\cite[page~4]{Ziegler1995}.  Suppose that $\tau
\subset \sigma$.  If $P(\tau)$ lies in the boundary of $P$ and if for
all $\tau \subsetneq \mu \subset \sigma$ we have $\dim P(\tau) < \dim
P(\mu)$ then we say $P(\tau)$ is a \textbf{face} of $P$.  Let
$P^{(k)}$ be the \textbf{ $k$--skeleton}: the union of the
$k$--dimensional faces of $P$.  A maximal chain of faces
\[
P(\tau_0) \subset P(\tau_1) \subset \ldots \subset P(\tau_n) = P
\]
is called a \textbf{flag}.  Then $P$ is \textbf{regular} if for any
two flags $F$ and $G$ of $P$ there is an isometry of $\R^n$ that
preserves $P$ and sends $F$ to $G$.

In dimensions one, two, and three the regular polytopes are known of
old.  These are the interval, the regular $k$--gons, and the
\textbf{Platonic solids}: the tetrahedron (simplex), the cube, the
octahedron (cross-polytope), the dodecahedron, and the icosahedron.
In all higher dimensions there are versions of the simplex, cube, and
cross-polytope.  In dimension four these are known as the $5$--cell,
the $8$--cell, and the $16$--cell.  Surprisingly, the only remaining
regular polytopes appear in dimension four!  There are only three of
them: the $24$--cell, the $120$--cell, and the
$600$--cell~\cite[page~136]{Coxeter1973}.


Suppose $P$ is a regular $n$--polytope.  The extreme symmetry of $P$
implies that we can move $P$ so that the vertices $P^{(0)}$ lie in the
unit sphere, $S^{n-1}$.  Projecting radially from the origin transfers
the one-skeleton $P^{(1)}$ from $\R^{n}$ into $S^{n-1}$.
Stereographic projection then places $P^{(1)}$ in $\R^{n-1}$.

Applied to a $4$--polytope, these projections turn the Euclidean
geometry of $P^{(1)}$ first into a design of arcs of great circles in
$S^3$ and then into a design of segments of circlines in $\R^3$.  If
$P^{(1)}$ meets the projection point then the design includes line
segments running off to infinity.  Coincidentally, Figure~\ref{axes.pdf} shows 
this for the $16$--cell.  In order to produce such a
design as a physical object, we need to thicken the circline segments
to have non-zero volume.
One possible approach uses the Euclidean geometry of $\R^3$: we could
thicken all segments of the design to get tubular neighbourhoods of
constant radius.  However, the result is not satisfactory; near the
origin in $\R^3$ the tubes are much too thick compared to their
separation.

\begin{figure}[htb]
\centering
\subfloat[A generic viewpoint.]{
\includegraphics[width=0.22\textwidth]{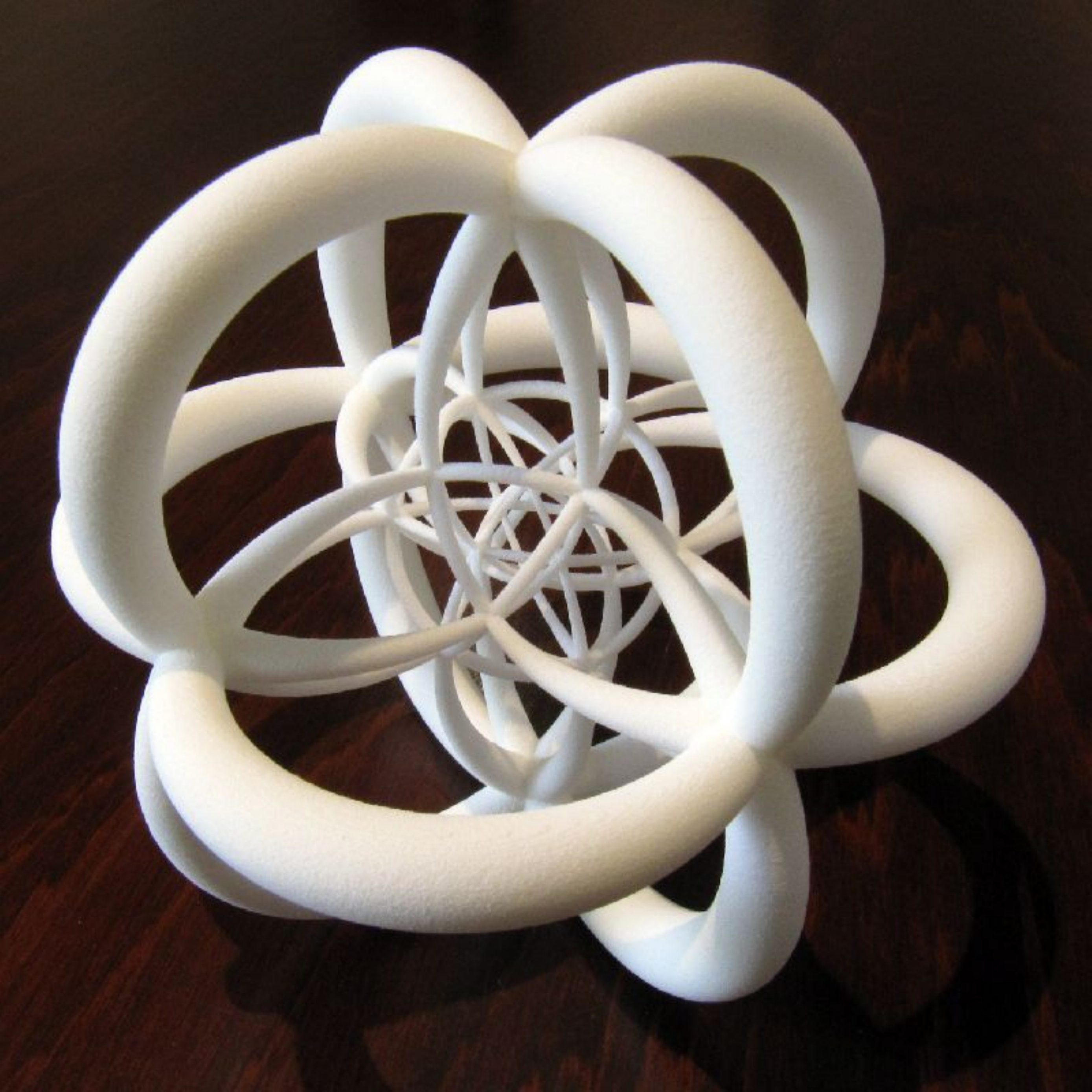}
\label{24_cell_1}}
\quad
\subfloat[A 2-fold symmetry axis.]{
\includegraphics[width=0.22\textwidth]{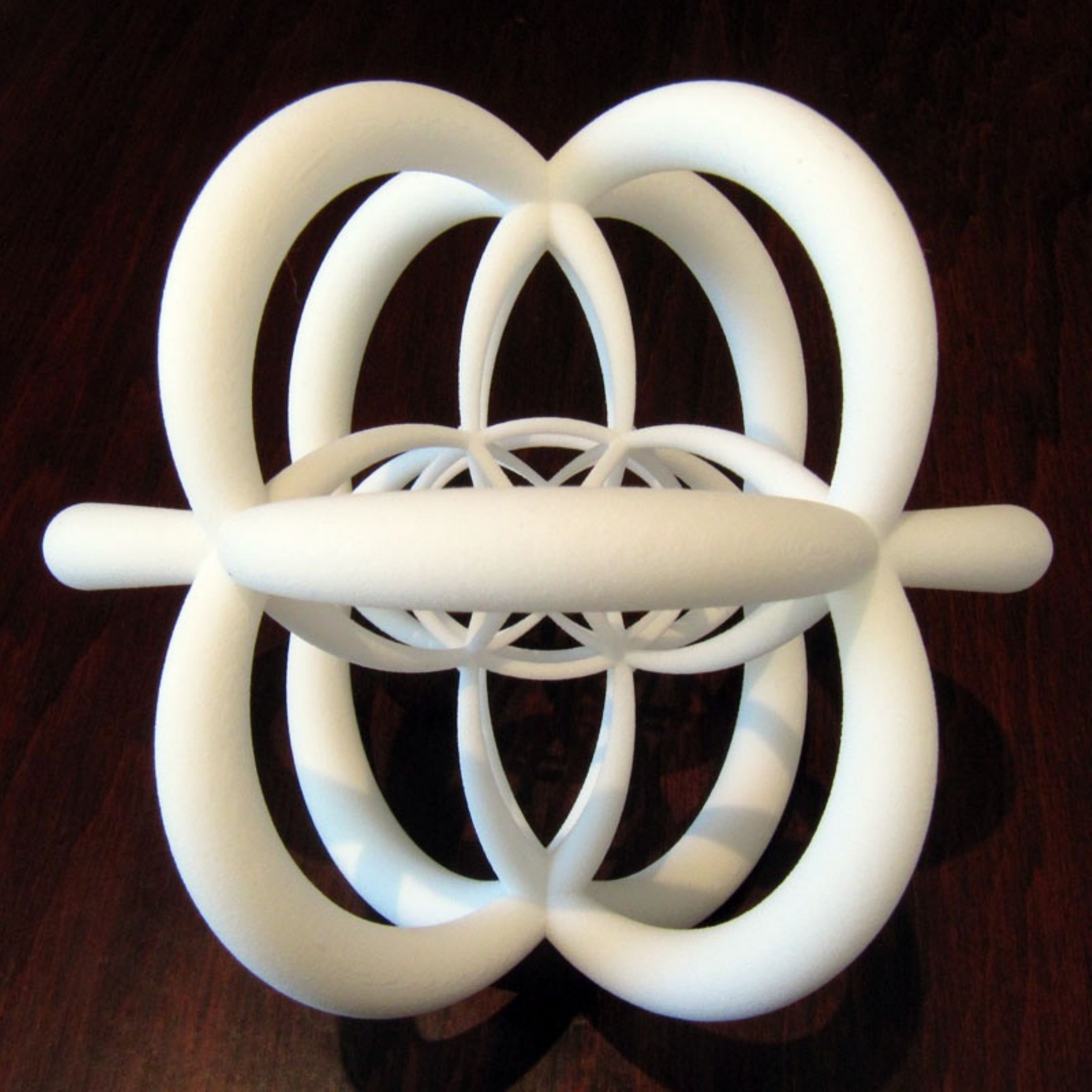}
\label{24_cell_symmetry2}}
\quad
\subfloat[A 3-fold symmetry axis.]{
\includegraphics[width=0.22\textwidth]{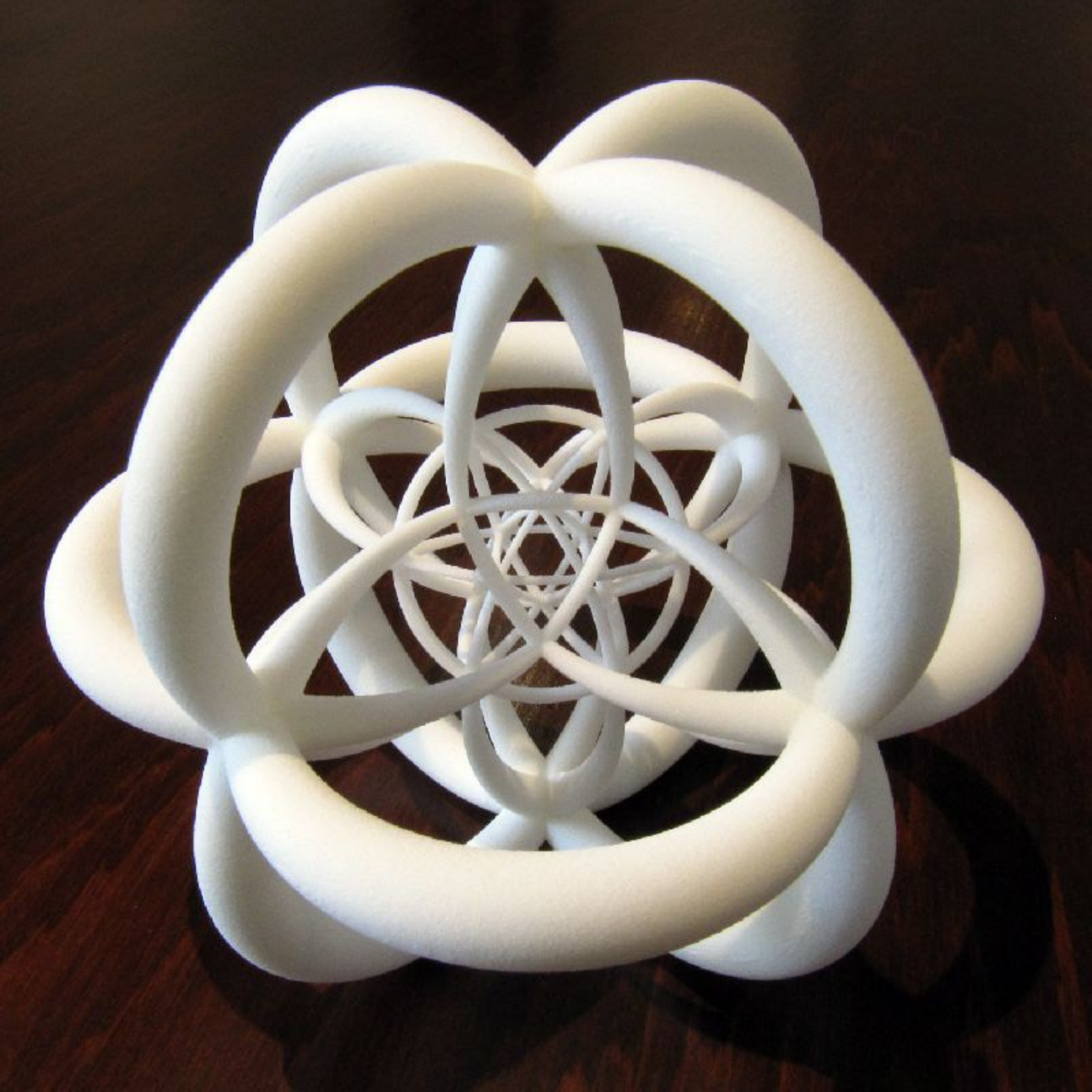}
\label{24_cell_symmetry3}}
\quad
\subfloat[A 4-fold symmetry axis.]{
\includegraphics[width=0.22\textwidth]{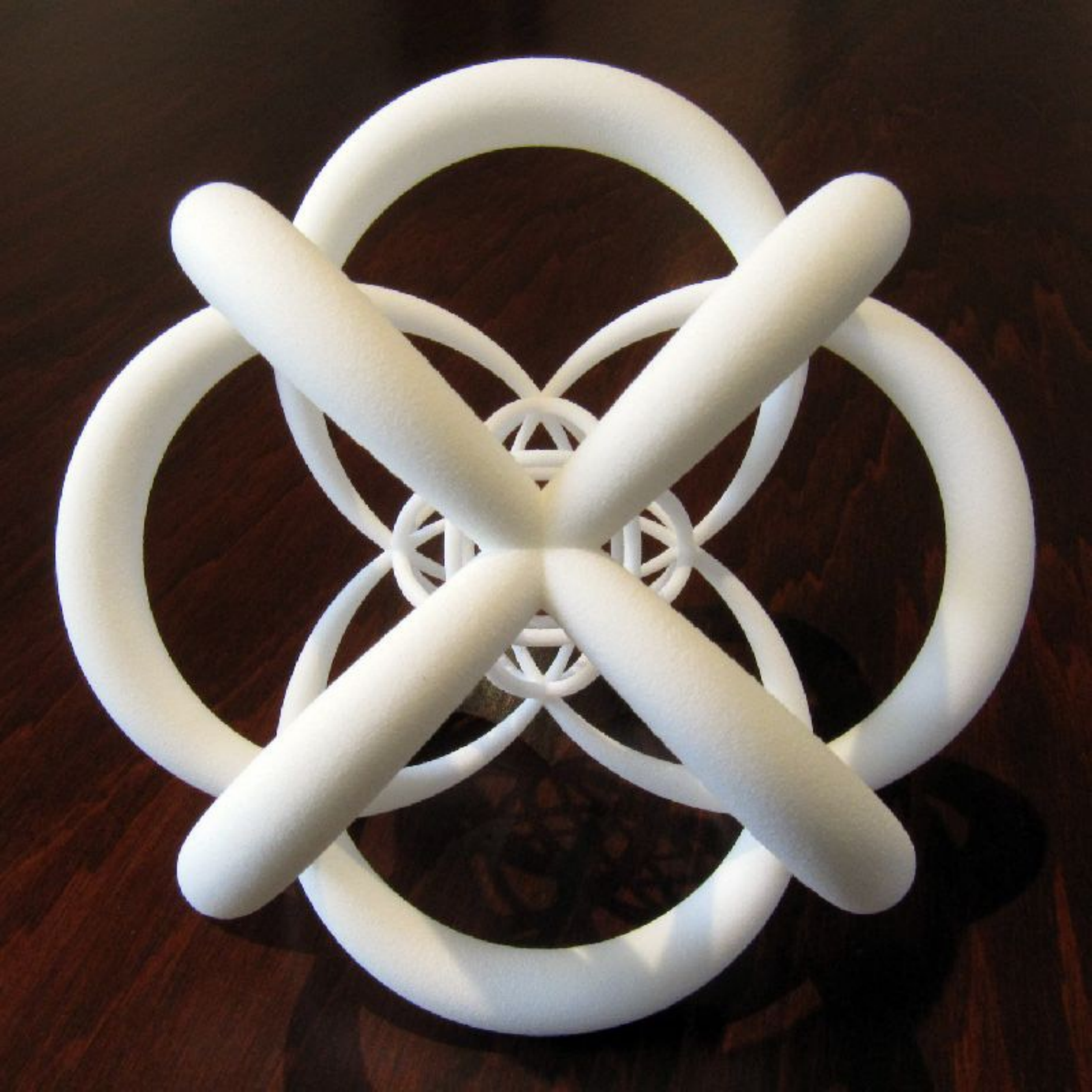}
\label{24_cell_symmetry4}}
\caption{$24$--Cell, 2011, $9.0\times 9.0 \times 9.0$ cm.}
\label{24-Cell}
\end{figure}

A better solution is to use tubular neighbourhoods in the intermediate
$S^3$ geometry.
For this we must parameterise the image of such a tube under
stereographic projection. Here the circline property is very
useful. The boundary of a tubular neighbourhood of a geodesic in $S^3$
can be made as a union of small circles in $\R^4$.  (These circles lie
in $S^3$, but are not great.)
The small circles map to circlines in $\R^3$, which can be directly
parameterised.  Computer visualisation of stereographic projections of
4-polytopes, in this style, are beautifully rendered by the program
Jenn3d~\cite{Jenn3d}.  In Figure~\ref{24-Cell} we show four views of a
3D print of the $24$--cell, with the projection point chosen to be at
the center of one of the cells. See also Ocneanu's ``Octacube''
\cite{octacube}.
 
%

The sculpture in Figure~\ref{24-Cell} illustrates a problem inherent
in 3D printing of stereographic projections.  Suppose that $P$ is a
symmetric design in $S^3$ and $Q = \rho(P)$ is the stereographic
projection from the north pole, $N$.  In this case the largest
features of $Q$ will correspond to the parts of $P$ closest to $N$.
These are the main contributers to volume and thus to cost.  The
smallest features of $Q$ will be roughly half the size of the parts of
$P$ nearest the south pole.
The 3D printing process places a lower bound on the size of the
smallest printable feature: current technology allows around $1$ mm.

Of course we may scale $Q$ in $\R^3$; scaling up ensures printability
while scaling down reduces volume.  Thus printability and cost are in
tension.  For example, if we rotate the $120$--cell so that $N$ lies
at the center of a dodecahedral face, and stereographically project,
then the largest feature is around $29.4$ times larger than the
smallest.  So here scaling to ensure printability also ensures
unaffordability.

One solution to this problem, as employed by Hart~\cite{Hart2011},
is use a \textbf{projective transformation} instead of stereographic
projection.  This takes a $4$--polytope to its \textbf{Schlegel
diagram}~\cite[page 133]{Ziegler1995}.  This is typically much more
compact. However, conformality is lost; the resulting figure distorts
both lengths and angles.

\begin{figure}[htb]
\centering
\subfloat[A generic viewpoint.]{
\includegraphics[width=0.22\textwidth]{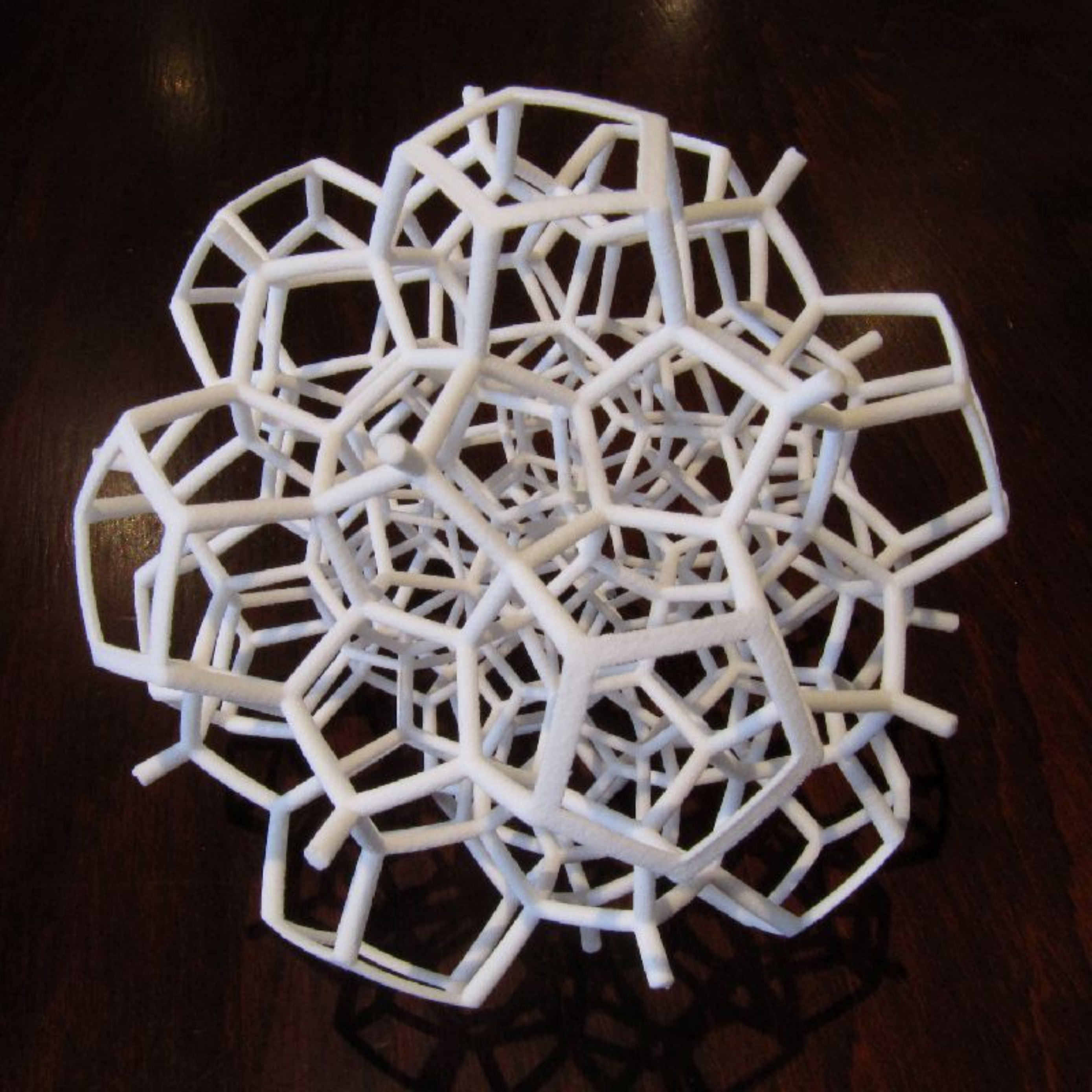}
\label{half_120-cell2}}
\quad
\subfloat[A 2-fold symmetry axis.]{
\includegraphics[width=0.22\textwidth]{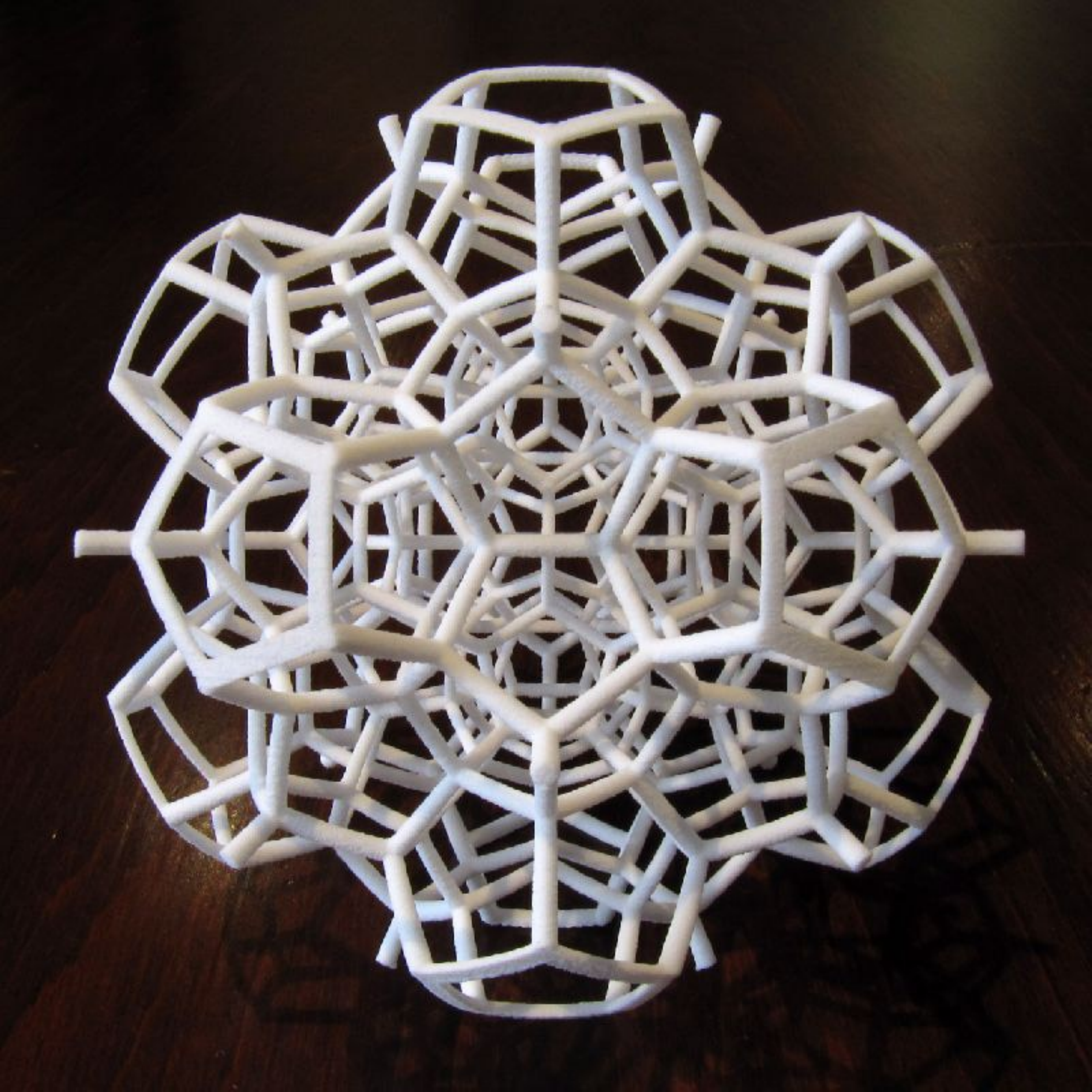}
\label{half_120-cell_sym_2}}
\quad
\subfloat[A 3-fold symmetry axis.]{
\includegraphics[width=0.22\textwidth]{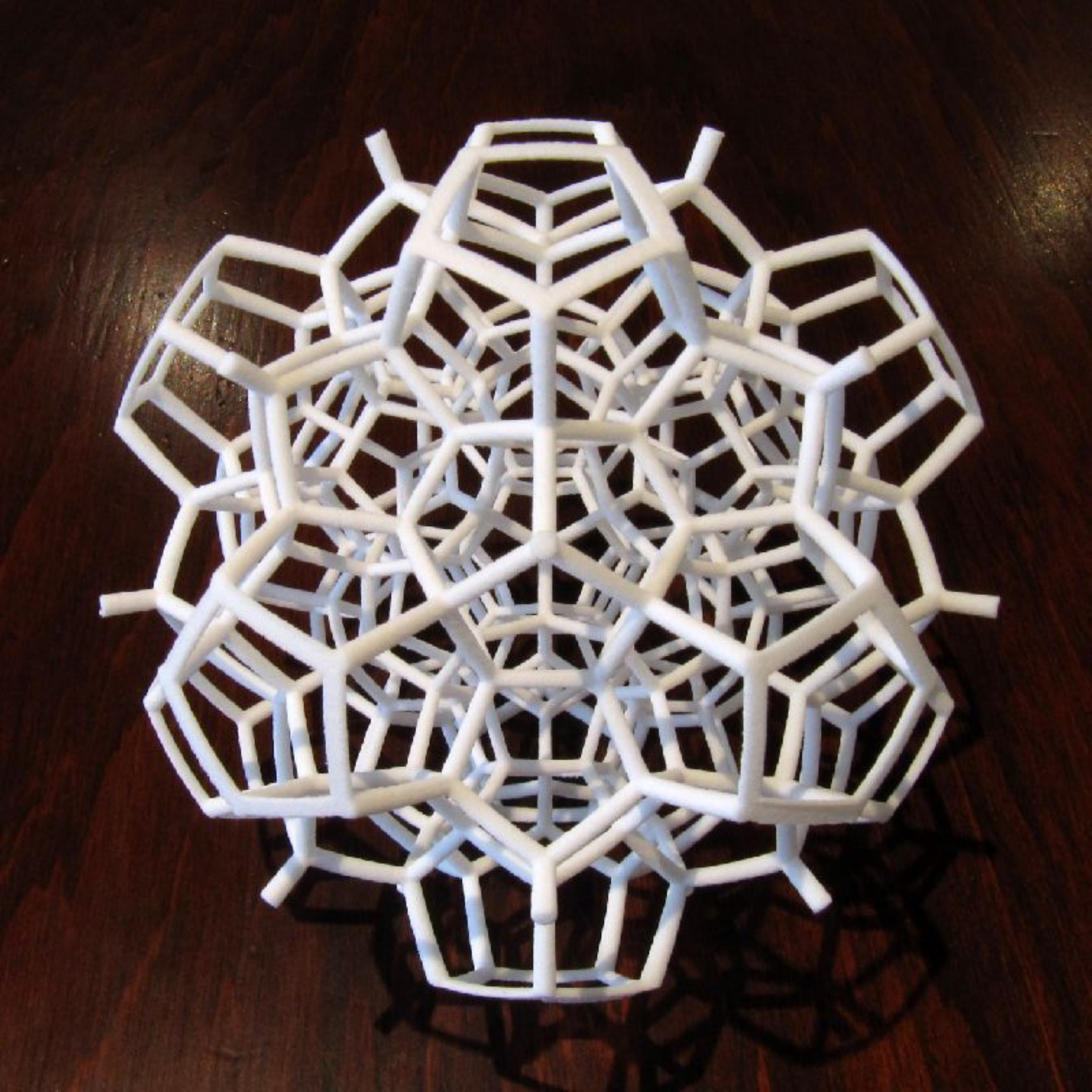}
\label{half_120-cell_sym_3}}
\quad
\subfloat[A 5-fold symmetry axis.]{
\includegraphics[width=0.22\textwidth]{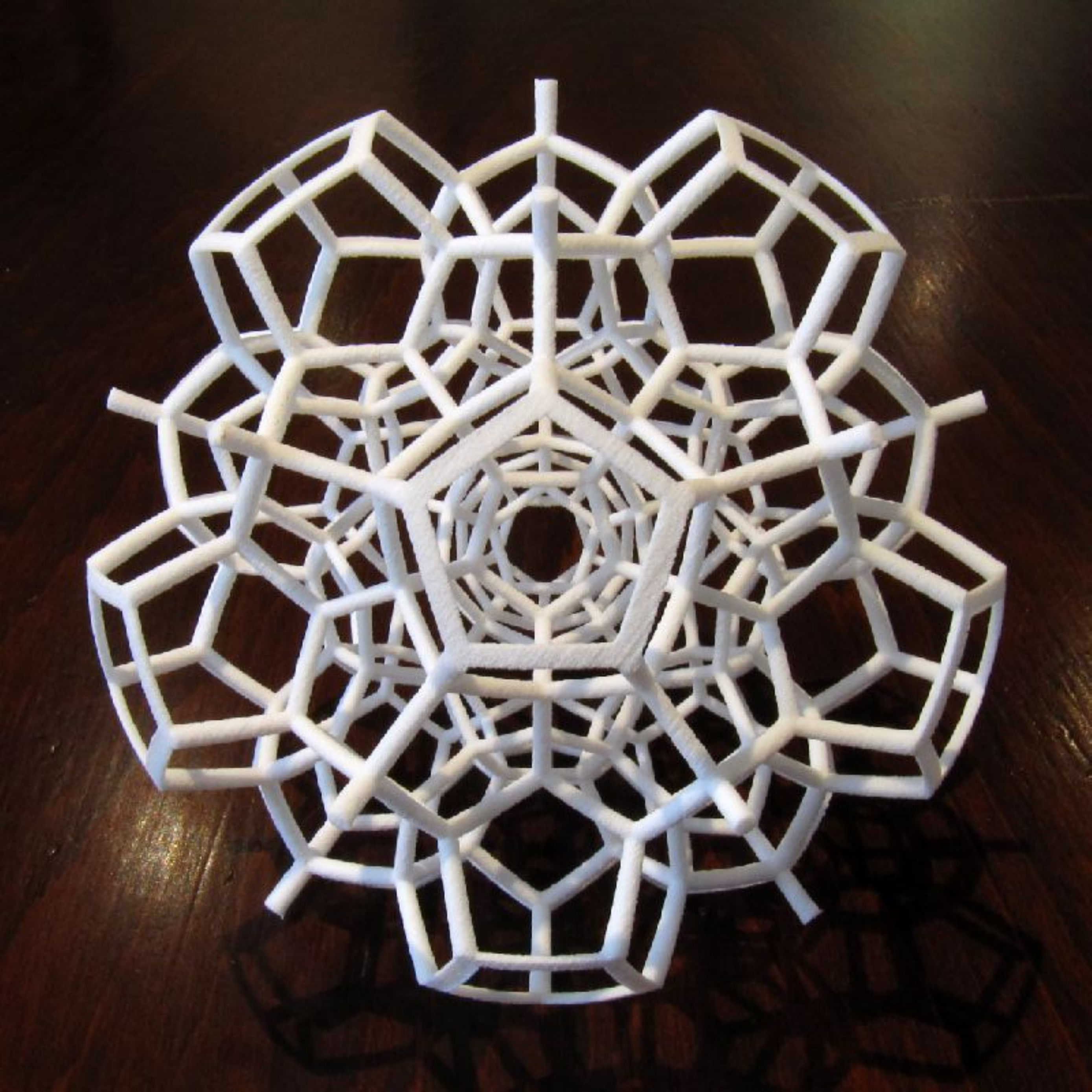}
\label{half_120-cell_sym_5}}
\caption{Half of a 120-Cell, 2011, $9.9\times 9.9 \times 9.9$ cm.}
\label{Half_120-Cell}
\end{figure}

Our alternative, shown in Figure~\ref{Half_120-Cell}, is to only print
half of the object.  We cut $S^3$ along the equatorial $S^2$; the
sphere equidistant from the north and south poles.  Choosing the north
pole as the projection point, we project the half of the design in the
southern hyperhemisphere.
The image is contained in the unit ball $\BB^3 = \{ x \in \R^3 : |x|
\leq 1 \}$. This done, the thinnest and thickest parts differ only by
a factor of two, at the most.  For stereographic projection, parts of
the design near the projection point are the real problem, in terms of
size.  Eliminating the half nearest the projection point eliminates
the problem.

%

\begin{figure}[htb]
\centering
\subfloat[A generic viewpoint.]{
\includegraphics[width=0.22\textwidth]{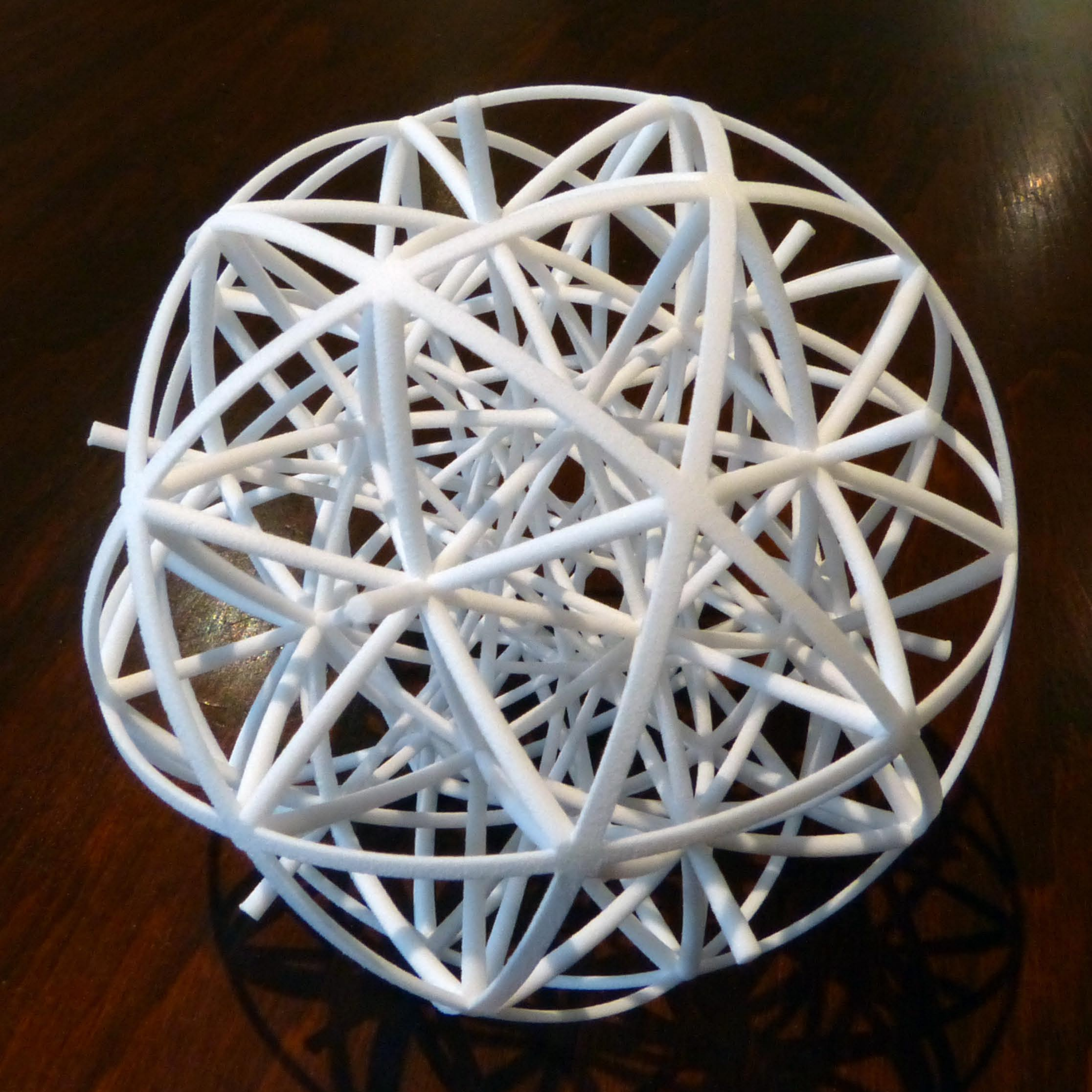}
\label{half_600-cell}}
\quad
\subfloat[A 2-fold symmetry axis.]{
\includegraphics[width=0.22\textwidth]{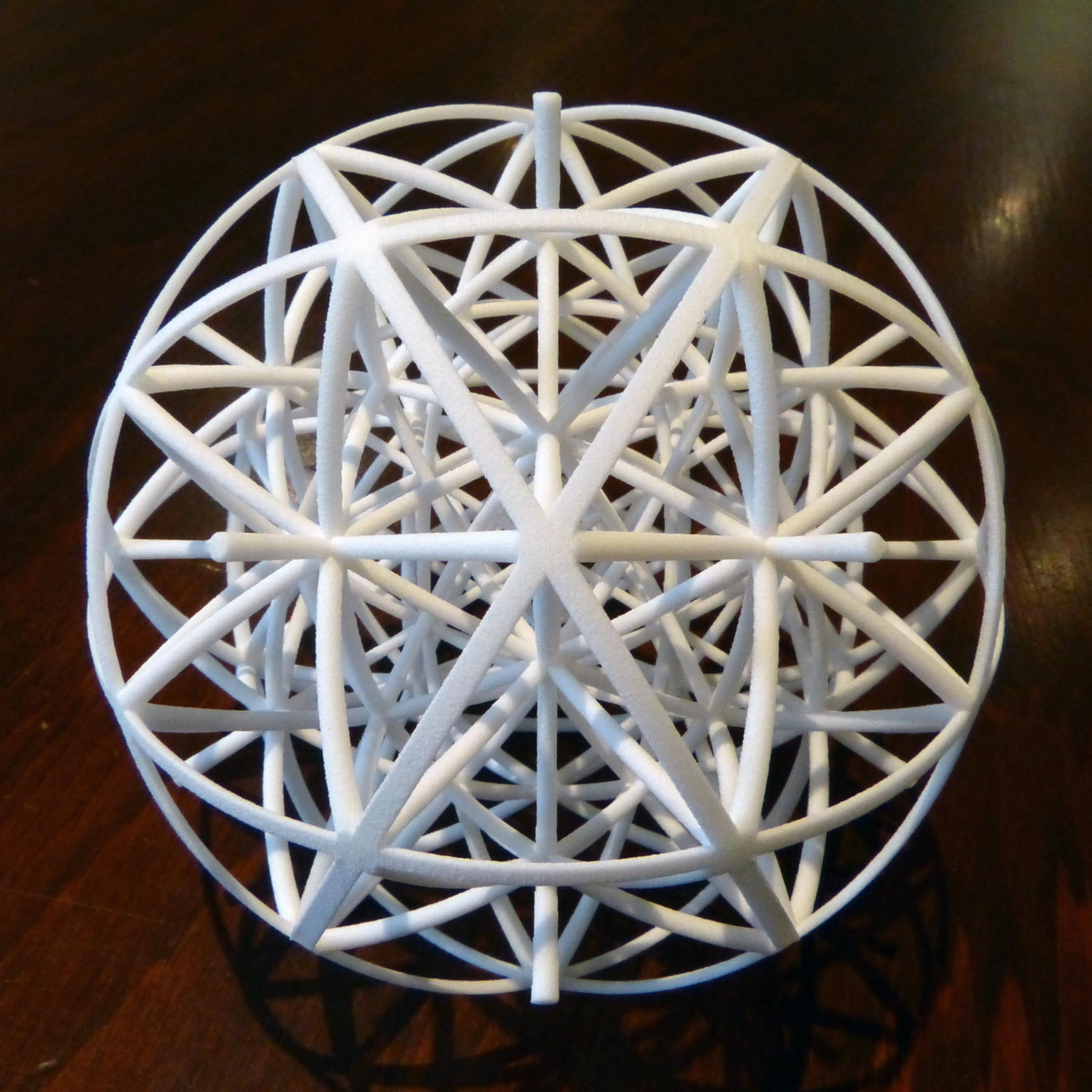}
\label{half_600-cell_sym_2}}
\quad
\subfloat[A 3-fold symmetry axis.]{
\includegraphics[width=0.22\textwidth]{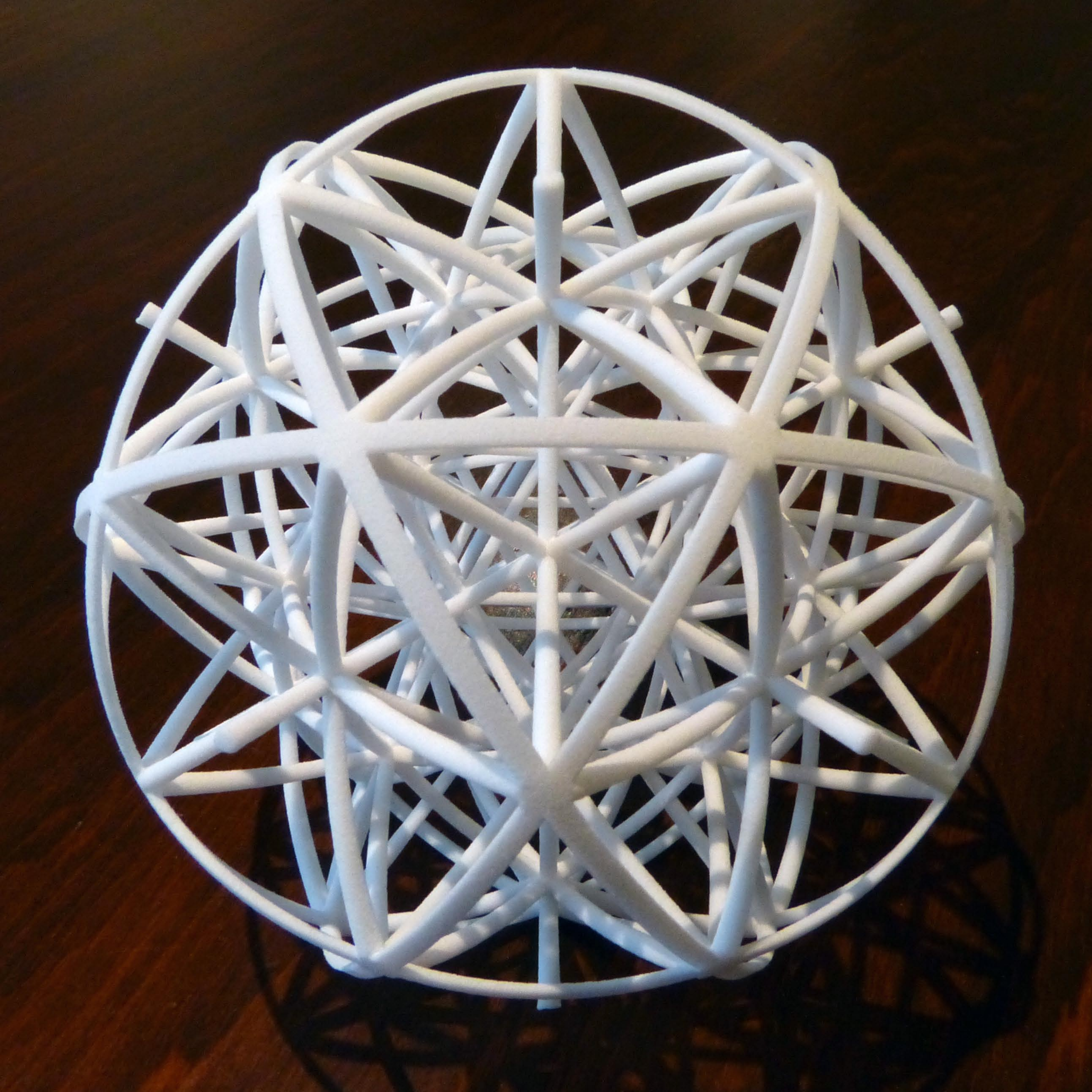}
\label{half_600-cell_sym_3}}
\quad
\subfloat[A 5-fold symmetry axis.]{
\includegraphics[width=0.22\textwidth]{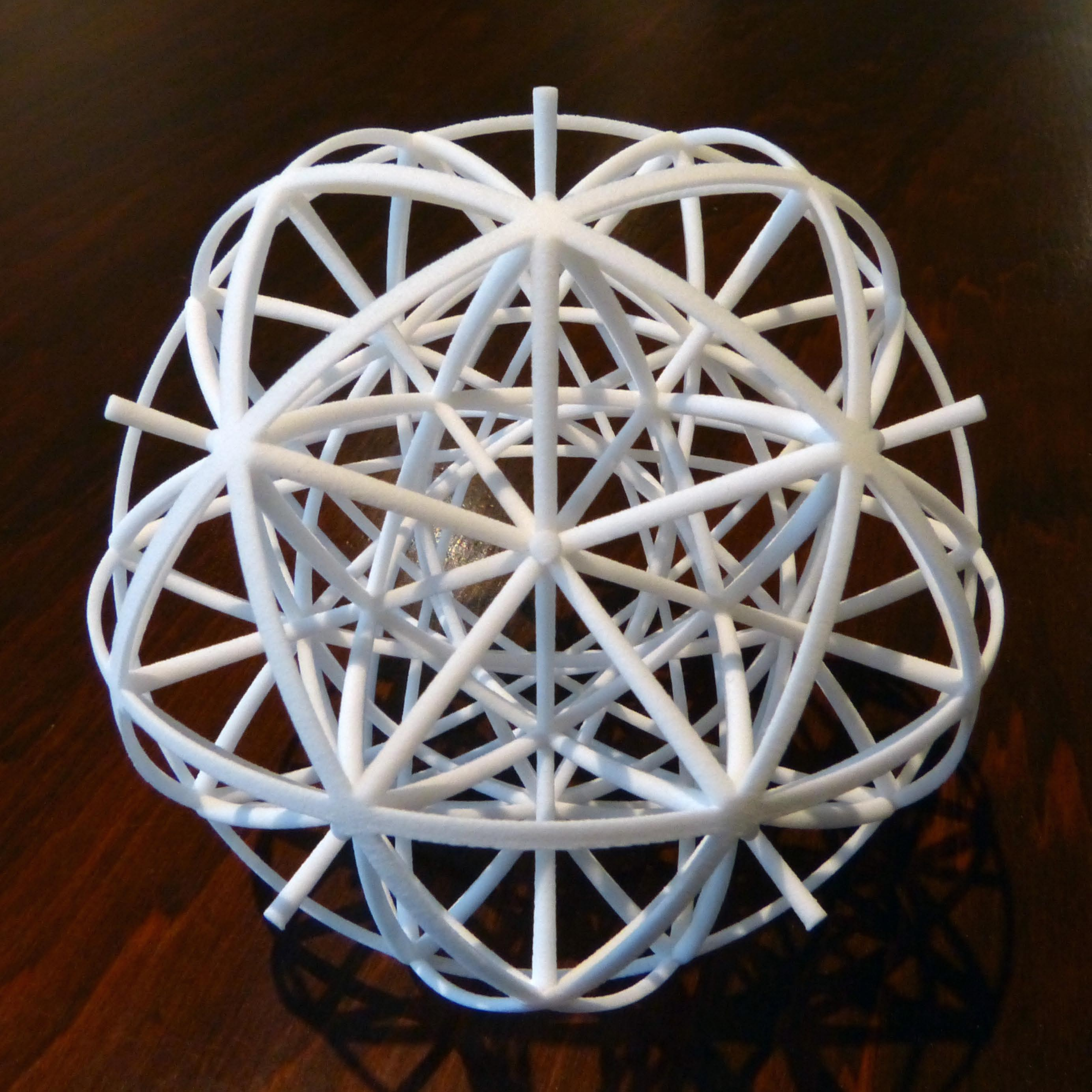}
\label{half_600-cell_sym_5}}
\caption{Half of a 600-Cell, 2011, $9.9\times 9.9 \times 9.9$ cm.}
\label{Half_600-Cell}
\end{figure}

Note that half of the $120$--cell is still very complicated!  However,
one can understand the whole of the $120$--cell by imagining
reflecting the object across the cutting two-sphere.  Note as well,
that printing only the southern half of a design allows us to print
objects that pass through the north pole, which ordinarily would be
infinitely expensive. For example, in Figure~\ref{Half_600-Cell} we
show one-half of the stereographic projection of the vertex centered
$600$--cell.  This version of the $600$--cell is positioned so as to
be dual to the facet-centered 120-cell shown in
Figure~\ref{Half_120-Cell}.  The other half of the vertex-centered
$600$--cell cannot be printed because the vertex antipodal to the
origin meets the projection point.

\subsection{Parameterisations of surfaces and torus knots}
\label{Sec:Tori}

The geometry of $S^3$ lends itself particularly well to the
representation of tori and torus knots.  There seem to be two reasons
for this.  First, in its natural position certain geodesics in the
torus are great circles in $S^3$.  Second, quaternionic multiplication
and its relatives directly parametrise torus knots.

When representing a surface as a 3D printed object, it is often a good
idea to drill holes in the surface, both to save on material used and
so the viewer can see, partly, through the surface to what is behind.
In our approach, the pattern of holes shows the parameterisation, by
realising the surface as a grid with grid-lines in the direction of
the parameters.

\paragraph{Clifford torus}
Recall that $e^{i\theta} = \cos(\theta) + i \sin(\theta)$ parametrises
a great circle $S^1$.  The same formula holds replacing $i$ everywhere
by $j$ or by $k$. 
A Clifford torus is foremost a torus, and so can be parameterised as a
product~\cite[page 139]{doCarmo1992} via
\begin{align*}
\TT = S^1 \times S^1 
& = 
\left\{
\frac{1}{\sqrt{2}}\bigl(
\cos(\alpha),\sin(\alpha),\cos(\beta),\sin(\beta)
\bigr)
\;\bigg\vert\;
0\leq \alpha <2\pi, \, 0\leq \beta < 2\pi
\right\}\\
& = 
\left\{
\frac{1}{\sqrt{2}}\bigl(e^{i\alpha} + e^{i\beta} \cdot j \bigr)
\;\bigg\vert\;
0\leq \alpha <2\pi, \, 0\leq \beta < 2\pi
\right\}.
\end{align*}
The factor of $1/\sqrt{2}$ rescales the torus to lie inside of the
unit sphere $S^3\subset\R^4$.  Note that if we vary $\alpha$ while
fixing $\beta$, then the point traces out a $(1,0)$ curve on $\TT$.
Conversely varying $\beta$ while fixing $\alpha$ yields a $(0,1)$
curve.  Unfortunately none of these curves are great circles in $S^3$.

On the other hand, if we vary $\alpha$ and $\beta$ simultaneously, at
the same (respectively, opposite) velocity the the point traces out a
$(1,1)$ (respectively $(1,-1)$) curve.  As we shall see, these are
great circles.

Note that $\TT$ divides $S^3$ into a pair of isometric \textbf{solid
tori}: copies of $S^1 \times D^2$.  We want to rotate the torus $\TT$
so that it meets the projection point.  This way, after stereographic
projection there is a pleasing symmetry; the two solid tori are
interchangeable by an isometry of $\R^3$.

We can use quaternions to fix the parameterisation, giving us the
desired $(1,1)$ and $(1,-1)$ curves, and to also move $\TT$ to meet
the projection point $1 \in S^3 \subset \H$.  Solving the second
problem first, note that $\frac{1}{\sqrt{2}}(1+j)$ lies in $\TT$.  If
$q$ is the quaternion satisfying $\frac{1}{\sqrt{2}}(1+j) q = 1$, then
$q = \frac{1}{\sqrt{2}}(1-j)$.  The new parameterisation of the torus
is given by post-multiplication by $q$:
\[
\frac{1}{\sqrt{2}}\bigl(e^{i\alpha} + e^{i\beta} \cdot j \bigr)
       \cdot \frac{1}{\sqrt{2}}(1-j)
= \frac{1}{2} \bigl( e^{i\alpha} + e^{i\beta} 
   + (e^{i\beta} - e^{i\alpha}) \cdot j \bigr).
\]
The torus meets the desired projection point when $\alpha = \beta =
0$.

We now solve the second problem, by rotating the coordinates through
$45^\circ$.  Take new coordinates $\theta, \phi$ where $\theta =
(\alpha+\beta)/2$ and $\phi = (\alpha-\beta)/2$.  So $\alpha = \theta
+ \phi$ and $\beta = \theta - \phi$.  Plugging in and simplifying, the
above parametrisation becomes $e^{i\theta} e^{-k\phi}$.  Keeping
$\phi$ fixed and varying $\theta$ now gives a $(1,1)$ curve, which is
also a great circle.  Note that we only need $0 \leq \theta < 2\pi, \,
0 \leq \phi < \pi$ to cover the whole torus.
We permute coordinates and change a sign to get a slightly neater
form:
\[
e^{i\phi}e^{j\theta} = 
\bigl(\cos(\theta)\cos(\phi), \cos(\theta)\sin(\phi), 
      \sin(\theta)\cos(\phi), \sin(\theta)\sin(\phi)\bigr)
\]
for $0 \leq \theta < 2\pi, 0 \leq \phi < \pi$.  The operations of
permuting the coordinates and changing the sign are symmetries of
$S^3$, so the geometry is unchanged and the surface $\TT$ still meets
the desired projection point, $1$.  The resulting parametrization is
Lawson's minimal surface $\tau_{1,1}$; see~\cite{lawson}.


\begin{figure}[htb]
\centering
\subfloat[A 2-fold symmetry axis.]{
\includegraphics[width=0.428\textwidth]{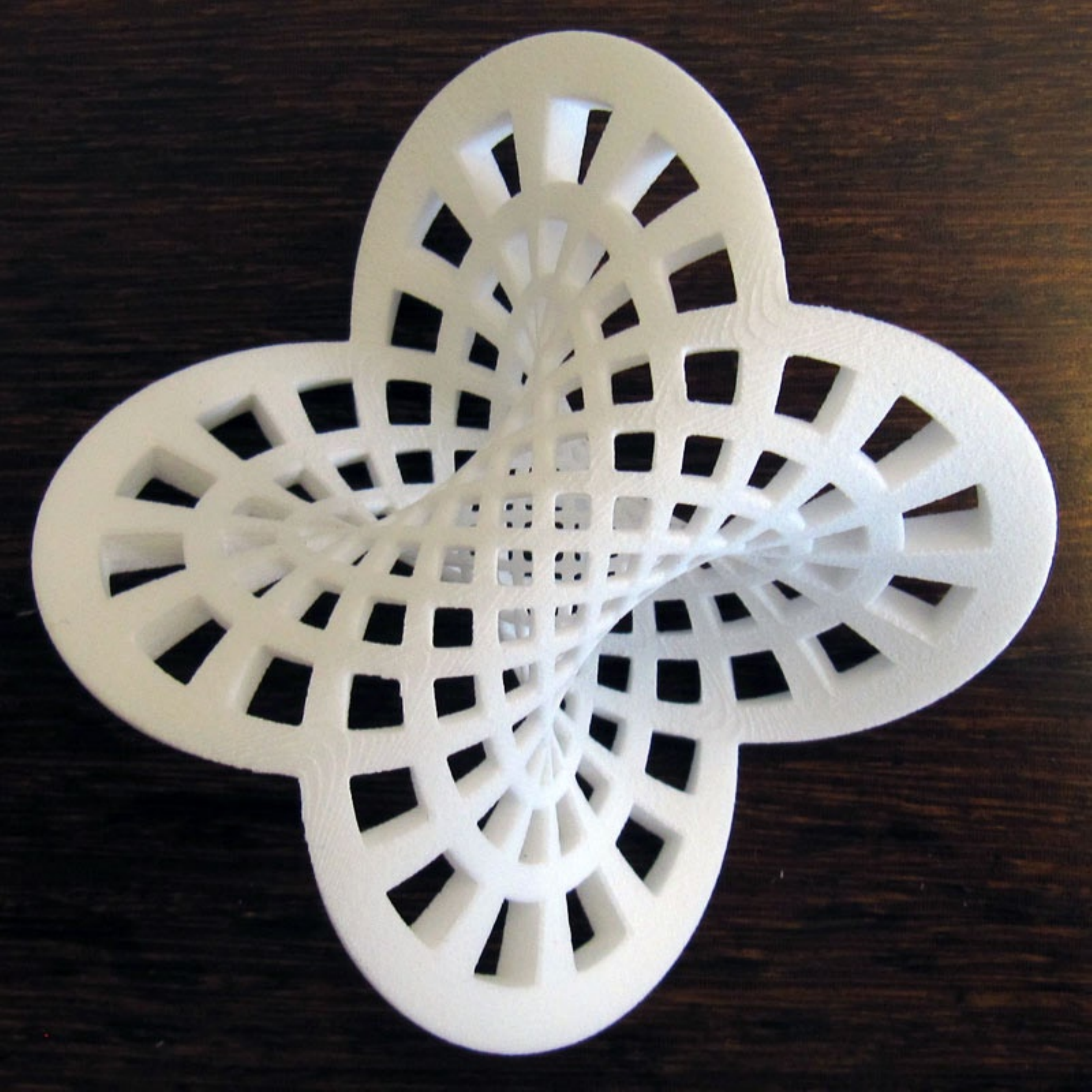}
\label{hopf_fibration_grid2_small}}
\qquad
\subfloat[A generic viewpoint.]{
\includegraphics[width=0.501\textwidth]{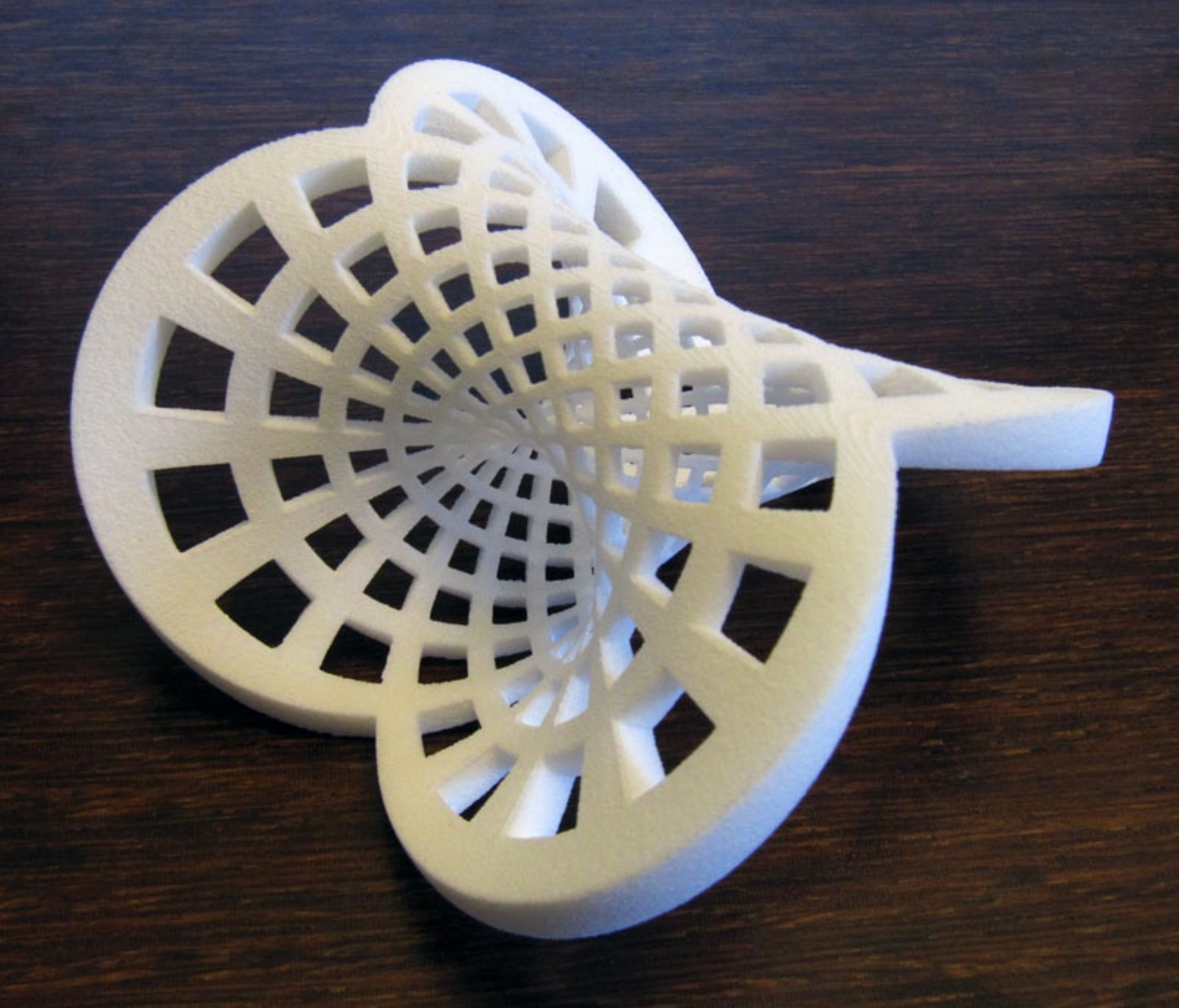}
\label{hopf_fibration_grid_small}}
\caption{Clifford Torus, 2011, $10.8\times 10.8 \times 3.4$ cm. }
\label{Hopf Fibration Grid}
\end{figure}

\paragraph{Finding the normal}
After stereographic projection, we get a 2-dimensional surface in
$\R^3 \cup \{\infty\}$.  As in Section~\ref{Sec:Polytopes}, for 3D
printing we must thicken our design to have positive volume.  Our plan
is to additionally parametrise the \textbf{normal} (that is,
perpendicular) to the surface, and then thicken in that
direction. As before, we do this thickening in $S^3$ rather than
$\R^3$.

Suppose that $F$ is any surface in $S^3$, with parametrisation
$p(\theta,\phi) \in S^3 \subset \R^4$.  Compute the tangent vectors
$\frac{\partial}{\partial\theta}p(\theta,\phi)$ and
$\frac{\partial}{\partial\phi}p(\theta,\phi)$ in $\R^4$.  Since $F$ lies
in $S^3$, these vectors are tangent to $S^3$ and so perpendicular to
$p(\theta,\phi)$, thought of as a vector from the origin.  The desired
normal vector $n(\theta, \phi)$ is a unit vector that is perpendicular
to the three given vectors $p$, $\frac{\partial}{\partial\theta}p$,
and $\frac{\partial}{\partial\phi}p$.  This determines $n$ up to sign.
Thus finding $n$ amounts to computing the kernel of the matrix with
rows $p$, $\frac{\partial}{\partial\theta}p$ and
$\frac{\partial}{\partial\phi}p$.  As these vectors vary with the
parameters $\theta$ and $\phi$ it is most convenient to compute $n$
via an application of Cramer's rule: $n$ is the determinant of the
matrix with first three rows $p$, $\frac{\partial}{\partial\theta}p$,
$\frac{\partial}{\partial\phi}p$, and fourth row the vector $(1, i, j,
k)$.

For the above parametrisation of the Clifford torus we find:
\[
\begin{array}{rccrrrrc}
p(\theta,\phi) & 
   = & \bigl(& \cos(\theta)\cos(\phi),&\cos(\theta)\sin(\phi),&\sin(\theta)\cos(\phi),&\sin(\theta)\sin(\phi)&\bigr)\\
\frac{\partial}{\partial\theta}p(\theta,\phi) & 
   = & \bigl(&-\sin(\theta)\cos(\phi),&-\sin(\theta)\sin(\phi),&\cos(\theta)\cos(\phi),&\cos(\theta)\sin(\phi)&\bigr)\\
\frac{\partial}{\partial\phi}p(\theta,\phi) & 
   = & \bigl(&-\cos(\theta)\sin(\phi),&\cos(\theta)\cos(\phi),&-\sin(\theta)\sin(\phi),&\sin(\theta)\cos(\phi)&\bigr)\\
n(\theta,\phi) & 
   = & \bigl(&-\sin(\theta)\sin(\phi),&\sin(\theta)\cos(\phi),&\cos(\theta)\sin(\phi),&-\cos(\theta)\cos(\phi)&\bigr)
\end{array}
\]


We introduce the parameter $\psi$ for the thickness of the
surface.  We move a distance $\psi$ along the geodesic from
$p(\theta, \phi)$ to $n(\theta, \phi)$ to reach
$r(\theta, \phi, \psi) = \cos(\psi)p(\theta,\phi) +
\sin(\psi)n(\theta,\phi)$.
Let $N_\epsilon(\TT)$ be the $\epsilon$--neighborhood of $\TT$, taken in $S^3$.  This is the same as thickening $\TT$ in the normal direction, using $r$.

Since $N_\epsilon(\TT)$ contains the projection point, the sculpture $\rho(N_\epsilon(\TT))$ would have infinite volume.
 We therefore remove a rectangular solid from $N_\epsilon(\TT)$; the boundary of the removed material is visible around the outside
 of the sculpture shown in Figure~\ref{Hopf Fibration Grid}.
 
\begin{figure}[htb]
\centering
\subfloat[A 2-fold symmetry axis.]{
\includegraphics[width=0.493\textwidth]{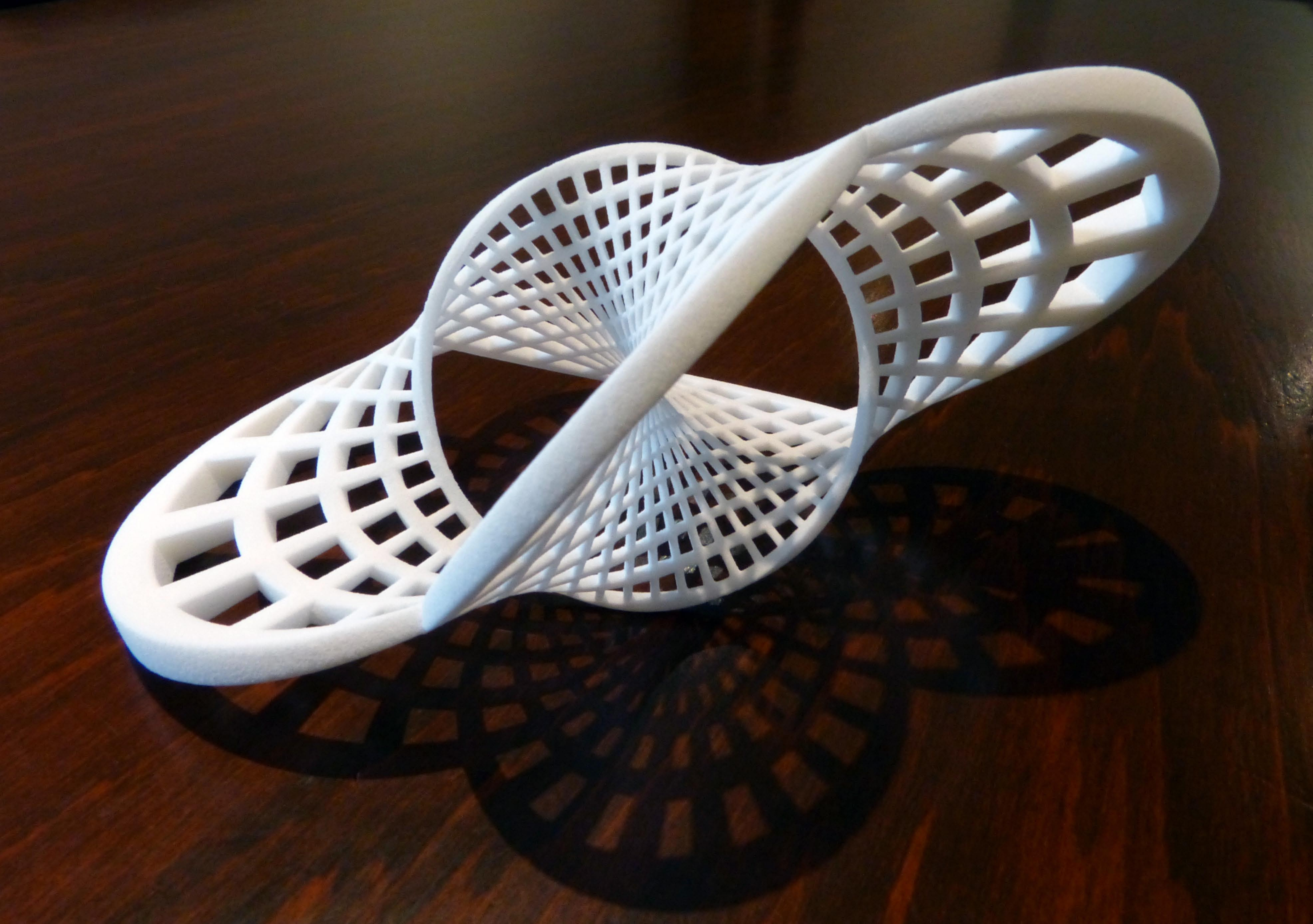}
\label{round_mobius_16_24_symmetry_small}}
\qquad
\subfloat[A generic viewpoint.]{
\includegraphics[width=0.437\textwidth]{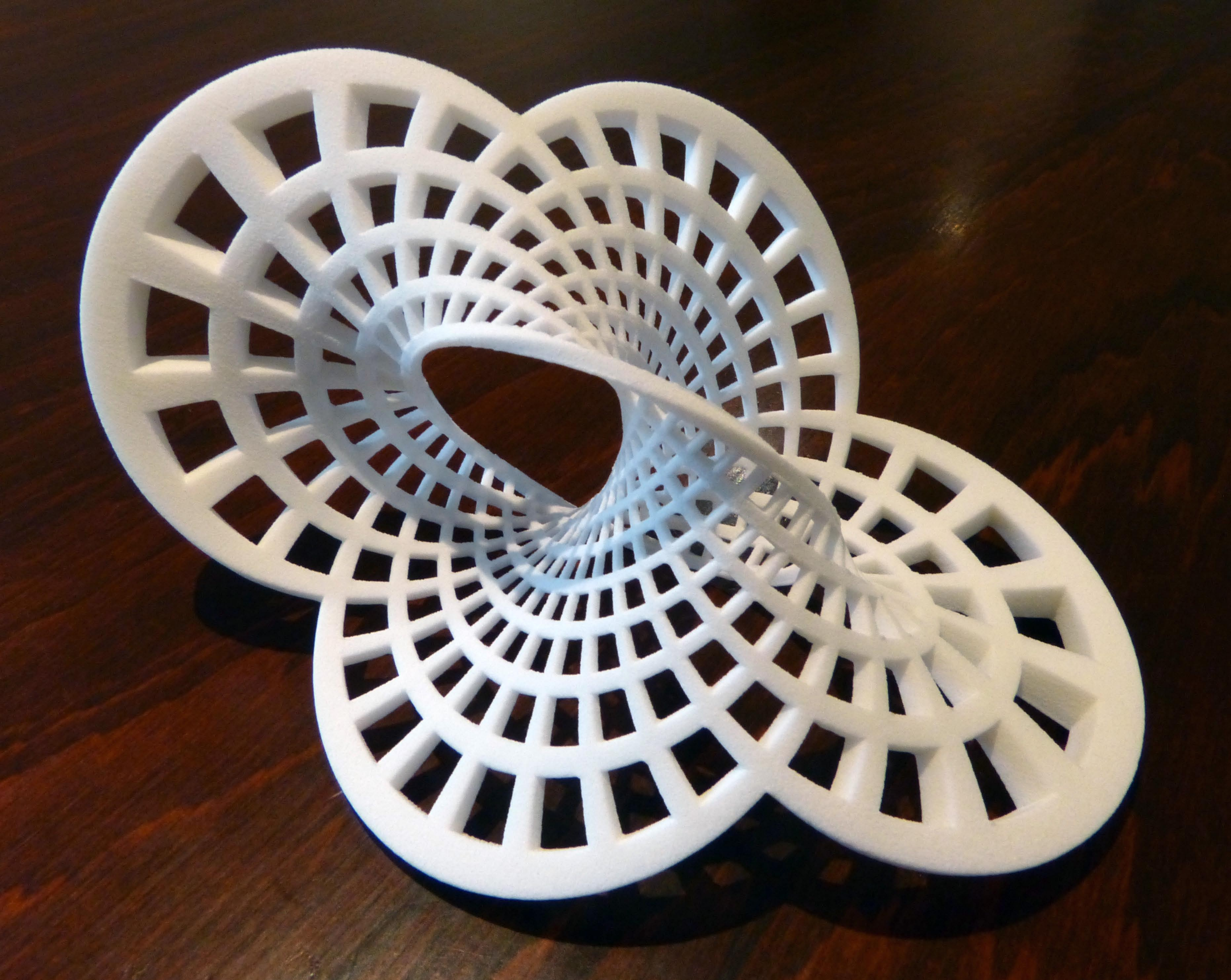}
\label{round_mobius_16_24_2_small}}
\caption{Round M\"{o}bius Strip, 2011, $15.2\times 10.9 \times 6.2$ cm.}
\label{Round M\"{o}bius Strip}
\end{figure}  
 
\begin{figure}[htb]
\centering
\subfloat[The 4-fold symmetry axis.]{
\includegraphics[width=0.242\textwidth]{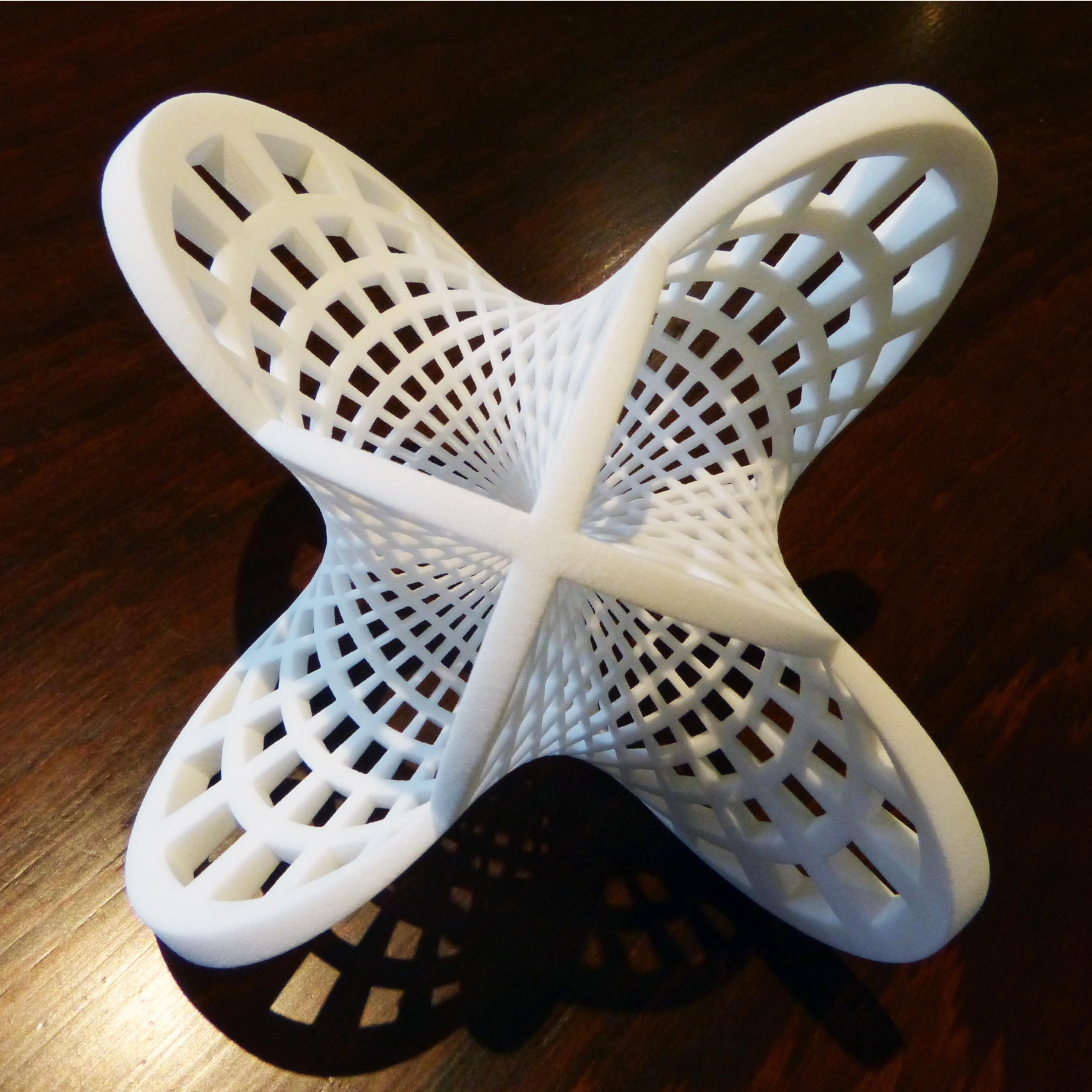}
\label{round_klein_16_24_symmetry4}}
\quad 
\subfloat[One of the 2-fold symmetry axes.]{
\includegraphics[width=0.347\textwidth]{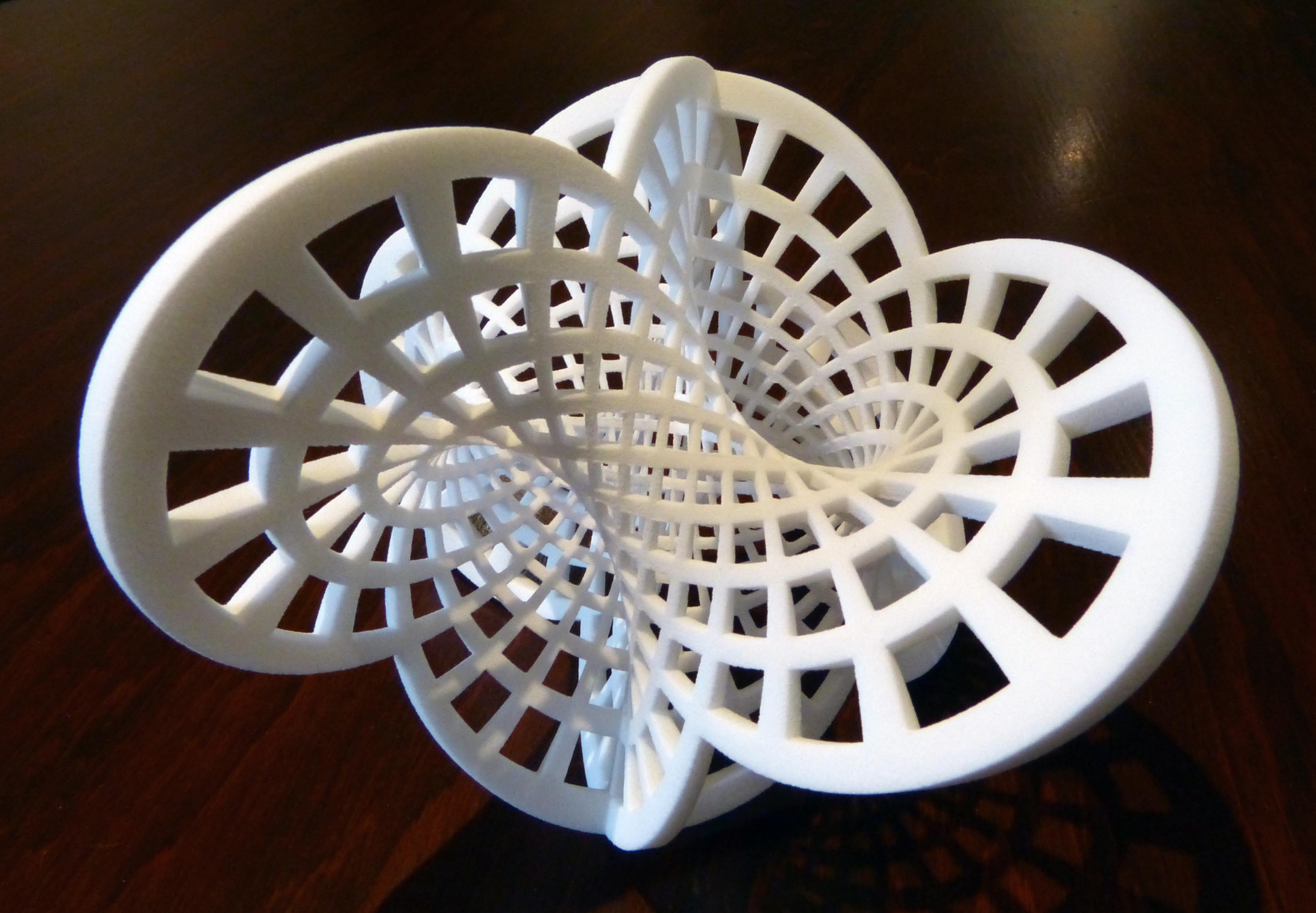}
\label{round_klein_16_24_symmetry2A}}
\quad
\subfloat[The other 2-fold symmetry axis.]{
\includegraphics[width=0.319\textwidth]{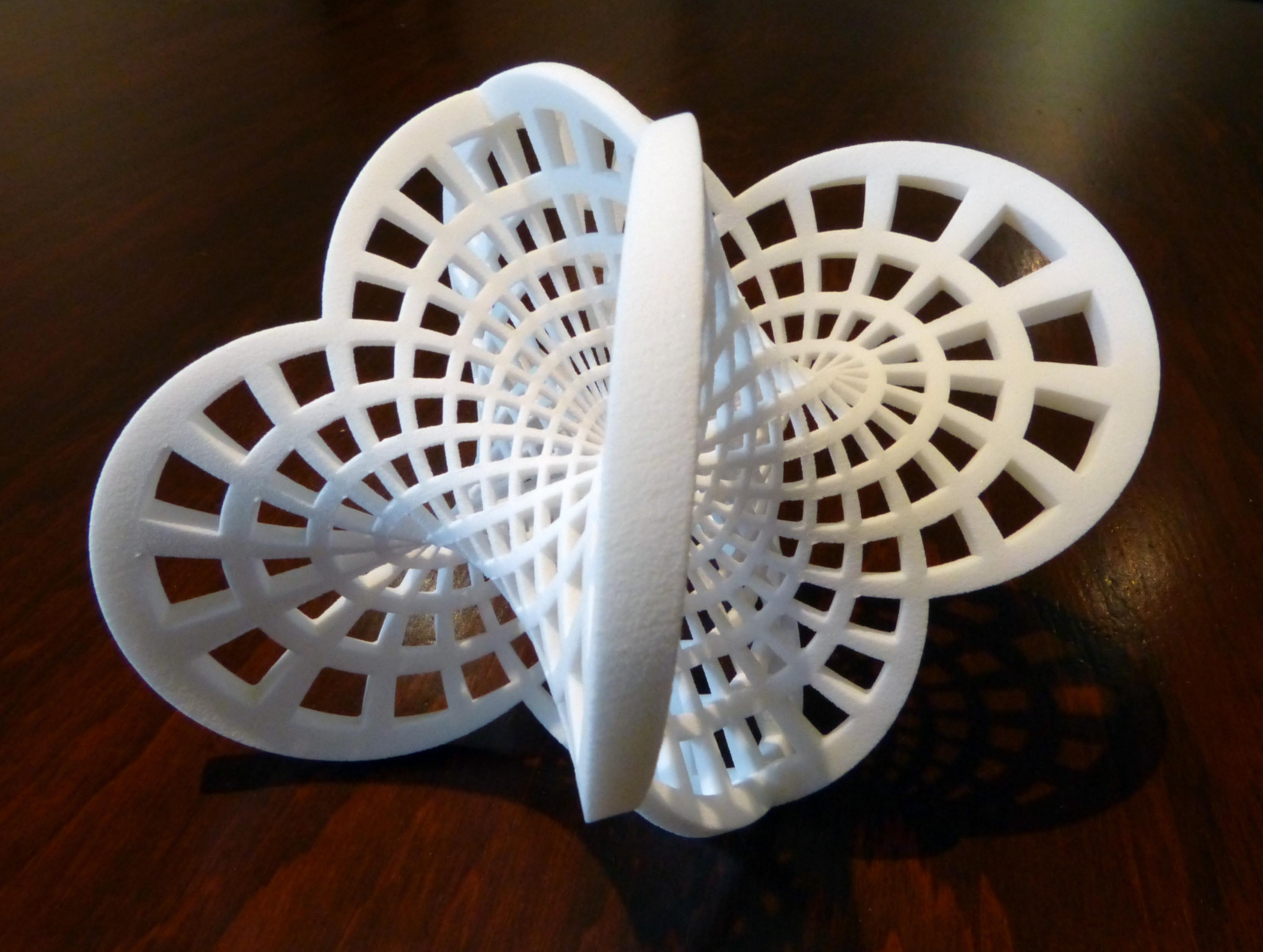}
\label{round_klein_16_24_symmetry2B}}
\caption{Round Klein Bottle, 2011, $15.2\times 15.2 \times 10.9$ cm.}
\label{Round Klein Bottle}
\end{figure} 
 
\paragraph{M\"{o}bius strip and Klein Bottle} 
A slight variant of the torus gives a M\"{o}bius strip:
\[
\left\{\bigl(
\cos(\theta)\cos(\phi),\cos(\theta)\sin(\phi),\sin(\theta)\cos(2\phi),\sin(\theta)\sin(2\phi)
\bigr) \mid 0\leq \theta <\pi,0\leq \phi < \pi\right\}
\]
This is a parameterisation of the ``Sudanese M\"{o}bius
strip''~\cite{sudanese}.  The border of the M\"{o}bius strip is given
by the points for which $\theta$ is $0$ or $\pi$.  Since these points
form a geodesic in $S^3$, the boundary is a circline in $\R^3$ by the
circline property. With the given parameterisation, stereographic
projection from $(0,0,-1,0)$ gives a circular boundary as in Figure \ref{Round M\"{o}bius
Strip}. The normal vector is calculated analogously to the torus case,
as follows. 
\[
\begin{array}{rcrrrrrc}
p(\theta,\phi) & 
  = & \bigl(& \cos(\theta)\cos(\phi),&\cos(\theta)\sin(\phi),&\sin(\theta)\cos(2\phi),&\sin(\theta)\sin(2\phi)&\bigr)\\
\frac{\partial}{\partial\theta}p(\theta,\phi) & 
  = & \bigl(&-\sin(\theta)\cos(\phi),&-\sin(\theta)\sin(\phi),&\cos(\theta)\cos(2\phi),&\cos(\theta)\sin(2\phi)&\bigr)\\
\frac{\partial}{\partial\phi}p(\theta,\phi) & 
  = & \bigl(&-\cos(\theta)\sin(\phi),&\cos(\theta)\cos(\phi),&-2\sin(\theta)\sin(2\phi),&2\sin(\theta)\cos(2\phi)&\bigr)\\
n(\theta,\phi) & 
  = & \frac{1}{\sqrt{1+3\sin^2(\theta)}}\bigl(&-2\sin(\theta)\sin(\phi),&2\sin(\theta)\cos(\phi),&\cos(\theta)\sin(2\phi),&-\cos(\theta)\cos(2\phi)&\bigr)
\end{array}
\]

Again the surface is punctured at the projection point, with a rectangular
hole in the grid pattern.  See Perry's sculpture ``Zero''~\cite{perry}
for a similar design. If we extend the strip across its boundary,
taking $0 \leq \theta < 2\pi$, we get the union of two punctured
M\"{o}bius strips, giving the twice-punctured Klein bottle shown in
Figure~\ref{Round Klein Bottle}.  This parameterisation of the
(unpunctured) Klein bottle is Lawson's surface $\tau_{2,1}$.

\paragraph{Torus knot} A further variant gives a parameterisation of a
torus knot, in this case the trefoil knot: 
\[
\left\{\bigl(
\cos(\theta)\cos(\phi),\cos(\theta)\sin(\phi),\sin(\theta)\cos((3/2)\phi),\sin(\theta)\sin((3/2)\phi)
\bigr) \mid 0\leq \phi < 4\pi\right\}
\]
Here $\theta$ has a fixed value, greater than 0 and smaller than
$\pi/2$. Altering the fraction $3/2$ will produce other torus
knots. The normal vector may be found as before; however for this
model we used an ``alternative'' to the normal vector, namely
\[
n(\theta,\phi)  =
\bigl(
-\sin(\theta)\sin(\phi),\sin(\theta)\cos(\phi),\cos(\theta)\sin((3/2)\phi),-\cos(\theta)\cos((3/2)\phi)
\bigr).
\]

\begin{wrapfigure}[13]{r}{0.32\textwidth}
\centering
\includegraphics[width=0.32\textwidth]{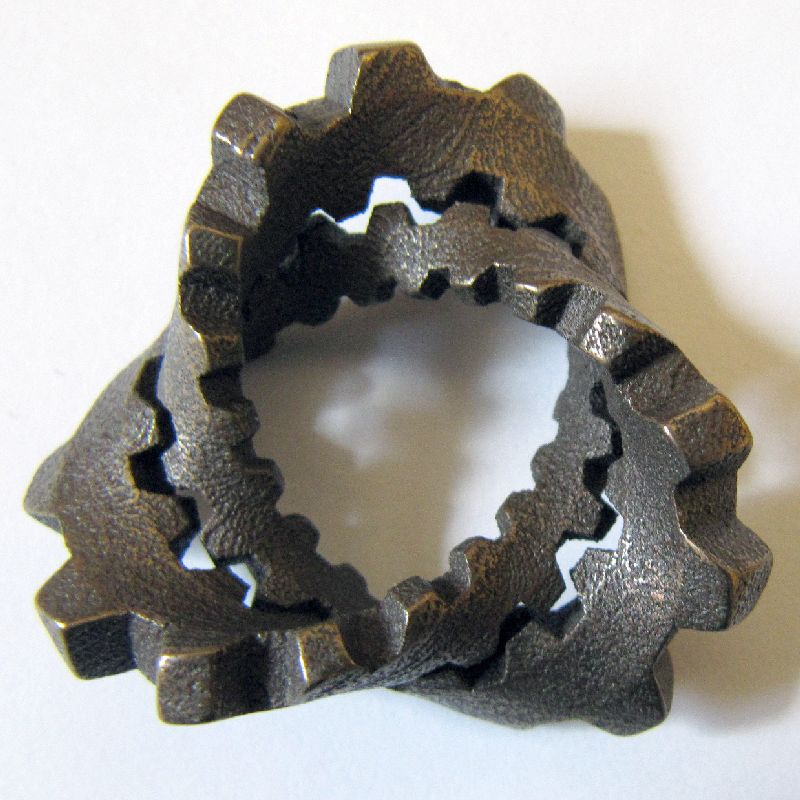}
\caption{Knotted Cog, 2011, $3.8\times 3.4 \times 1.3$ cm.}
\label{knotted_cog_stainless_steel_antique_bronze_glossy_square}
\end{wrapfigure}

Using the local coordinates $(\theta,\phi,\psi)$, we can add small
features to the sculpture, using any shape we could define in ordinary
3-dimensional space. In the case shown in
Figure~\ref{knotted_cog_stainless_steel_antique_bronze_glossy_square},
we add cog teeth, which are simply truncated pyramids in
$(\theta,\phi,\psi)$ coordinates. The alternative normal vector
adds a slight shear slope to the teeth, which we feel is aesthetically
preferable.

\section{Future directions}

Our sculptures are tangible representives of topological and geometric
abstractions.  In order to do this, we naturally must construct
designs that occur in $\R^3$: that is, in actual space.  In each case
we attempted to choose the most canonical such geometries available
and then the most faithful projections.

There is a wild array of further topological and combinatorial
objects.  For example, there is a rich theory of knots and surfaces
and their interrelations.  We have not yet found (or perhaps better,
understood) satisfactory geometric representations, or at least
representatives which map to $\R^3$ in satisfactory ways.  An example
of the latter problem would be surfaces of genus at least two.  These
have nice hyperbolic structures, but they cannot be mapped into $\R^3$
in a very satisfying way.


\setlength{\baselineskip}{13pt} 

\bibliographystyle{hamsplain}
\bibliography{spherebib}

\end{document}